\documentclass[a4paper,11pt]{amsart}

\usepackage{amssymb}
\usepackage{mhequ}
\usepackage[latin1]{inputenc}

\usepackage{cite}

\def \Sphere{\mathbb{S}}

\DeclareFontFamily{OT1}{rsfs}{}
\DeclareFontShape{OT1}{rsfs}{m}{n}{ <-7> rsfs5 <7-10> rsfs7 <10->
rsfs10}{} \DeclareMathAlphabet{\mathscr}{OT1}{rsfs}{m}{n}

\newcommand{\eq}[1]{\eqref{#1}}

\newcommand{\bel}[1]{\begin{equation}\label{#1}}
\newcommand{\bea}{\begin{eqnarray}}
\newcommand{\beaa}{\begin{eqnarray*}}
\newcommand{\bean}{\begin{eqnarray}\nonumber}
\newcommand{\beal}[1]{\begin{eqnarray}\label{#1}}
\newcommand{\beadl}[1]{\begin{deqarr}\label{#1}}
\newcommand{\eeadl}[1]{\arrlabel{#1}\end{deqarr}}
\newcommand{\eeal}[1]{\label{#1}\end{eqnarray}}
\newcommand{\eead}[1]{\end{deqarr}}
\newcommand{\eea}{\end{eqnarray}}
\newcommand{\eeaa}{\end{eqnarray*}}

\newcommand{\be}{\begin{equation}}
\newcommand{\ee}{\end{equation}}

\DeclareFontFamily{OT1}{rsfs}{}
\DeclareFontShape{OT1}{rsfs}{m}{n}{ <-7> rsfs5 <7-10> rsfs7 <10->
rsfs10}{} \DeclareMathAlphabet{\mycal}{OT1}{rsfs}{m}{n}

\newcounter{mnotecount}[section]

\newcommand{\rmnote}[1]{}

\def\mysavedown#1{\edef\mysubs{\mysubs#1}}
\def\mysaveup#1{\edef\mysups{\mysups#1}}
\def\mydown#1{{\mytensor}_{\vphantom{\mysubs}#1}}
\def\myup#1{{\mytensor}^{\vphantom{\mysups}#1}}
\def\tensor#1#2{
  #1
  \def\mytensor{\vphantom{#1}}
  \def\mysubs{\relax}
  \def\mysups{\relax}
  \let\down=\mysavedown
  \let\up=\mysaveup
  #2
  \let\down=\mydown
  \let\up=\myup
  #2
  }

\newcommand{\Tr}{\operatorname{Tr}}

\newcommand{\R}{\mathbb R}

\renewcommand{\div}{\operatorname{div}}
\renewcommand{\epsilon}{\varepsilon}
\renewcommand{\hat}{\widehat}

\def\crn#1#2{{\vcenter{\vbox{
        \hbox{\kern#2pt \vrule width.#2pt height#1pt
           }
          \hrule height.#2pt}}}}

\newcommand{\U}{\mathbb U}

\newcommand{\m}{\mathcal M}
\newcommand{\cH}{\mathcal H}

\newcommand{\w}{\widetilde}

\newcommand{\pre}[2]{{{\vphantom{#2}}^{#1}}\kern-.2ex{#2}}

\theoremstyle{plain}
\newtheorem{theorem}{\sc Theorem}[section]

\newtheorem{lemma}[theorem] {\sc Lemma}

\newtheorem{proposition}[theorem]{\sc Proposition}

\newtheorem{corollary}[theorem] {\sc Corollary}

\theoremstyle{definition}

\newtheorem{remark}[theorem]{\sc  Remark\rm}

\numberwithin{equation}{section}

\date{February 6, 2008}

\begin{document}

\title[Relative radial mass and Rigidity of some Warped Product Manifolds]
{Relative radial mass and rigidity of some Warped Product Manifolds }
\author[M. Arcostanzo]{Marc Arcostanzo}\address{Marc Arcostanzo,
Laboratoire d'analyse non lin\'eaire et g\'eom\'etrie, Facult\'e
des Sciences, 33 rue Louis Pasteur, 84000 Avignon, France}
\email{Marc.Arcostanzo@univ-avignon.fr}
\author[E. Delay]{Erwann Delay} \address{Erwann Delay, Laboratoire d'analyse
non lin\'eaire et g\'eom\'etrie, Facult\'e des Sciences, 33 rue
Louis Pasteur, 84000 Avignon, France}
\email{Erwann.Delay@univ-avignon.fr}
\urladdr{http://www.math.univ-avignon.fr/Delay}
\begin{abstract}
We give a Riccati type formula adapted for two metrics having
the same geodesics rays starting from a point or orthogonal to
an hypersurface, one of these metrics being a warped product
if the dimension $n$ is greater than or equal to 3.
This formula has non-trivial geometric consequences such as
a positive mass type theorem and other rigidity results. We
also apply our result to some standard models.
\end{abstract}

\maketitle

\\
\newline

\tableofcontents
\section{Introduction}\label{section:intro}

Let $n\geq2$ be an integer and
let $M=(a,b)\times N$, where $N$ is a compact $(n-1)$-dimensional manifold
(with or without boundary) and $(a,b)$ an open interval of $\R$.
We shall assume for simplicity that all objects defined on $M$ are of class
$C^{\infty}$, although this is far from necessary.
We consider a Riemannian  metric of the form
$$
g=dr^2+G(r)
$$
on $M$, where $G(r)$ is a family of Riemannian metrics on $N$.
We also consider a background warped product metric on M:
$$
g_0=dr^2+\mathring{G}(r),
$$
where $\mathring{G}(r)=\mathring{h}^2(r)\hat{G}$, $\mathring{h}$ being a positive
function on  $(a,b)$ and $\hat{G}$ a fixed metric on $N$.
For $(r,x)\in  M$, we define the  radial mass of $g$, relative to $g_0$,  at $(r,x)$ as
\footnote{In fact, only  $\mathring{h}$ appear in the dependance on $m$  upon $g_0$
but we choose this terminology for convenience}
\bel{masse}
m(r,x)=\mathring{h}^{2}(\mathring{H}-H),
\ee
where $H$ and $\mathring H$ are the mean curvatures of the level
set $N(r)$ of $r$ relative to $g$ and $g_0$ respectively.
We can already remark (see remark \ref{remdet}) that
\bel{masse2}
m(r,x)=
\mathring{h}^{2}\partial_r\left[\ln\left( \sqrt{ \frac{\det(g_0)}{\det(g)} } \right)\right].
\ee

We will first show the following Riccati type formula
\bel{riccadni}
\partial_rm=\mathring{h}^{2}(R_{rr}-\mathring{R}_{rr})+\frac{1}{n-1}\frac{m^2}{\mathring{h}^2}
+\mathring{h}^{2}|s|^2,
\ee
where $s$ is the trace free part of the $g-$Weingarten map of $N(r)$, whereas
$R_{rr}$ and $$\mathring{R}_{rr}=-(n-1)\frac{\mathring{h}''}{\mathring{h}}$$ are
the Ricci curvatures in the $(\partial_r,\partial_r)$ direction of $g$ and $g_0$
respectively.

This formula has many geometric consequences, the most interesting
being the following one. Let us define the radial mass of $g$, relative to $g_0$,
on $N(r)$ as

\bel{massint}
{\mathcal M} (r)=\int_Nm(r,x)d\mu_{\hat{G}}(x).
\ee

We will say that the metric $g$ is asymptotic to $g_0$ at $r_0\in[a,b]$ if near $r_0$,
$$
||G-\mathring{G}||_{\mathring{G}}=o(1),
$$
where $o(1)$ is uniform on $N$ (if $r_0 \in ]a,b[$, this means that the two metrics
coincide on $N(r_0)$).
Our main result can now be stated.
\begin{theorem}\label{thintro} If the limits of $\m$ exist at $a$ and $b$ and
$\mathring{h}^{2}(R_{rr}-\mathring{R}_{rr})\in L^1(drd\mu_{\hat{G}})$  then
$$
\m(b)-\m(a)\geq \int_M\mathring{h}^{2}(R_{rr}-\mathring{R}_{rr})drd\mu_{\hat{G}}$$
with equality if and only if $g$ is a warped product (i.e. $G=\mathring{h}^2\w{G}$,
$\w{G}$ being a fixed metric on $N$). Let us furthermore assume that $g$ is asymptotic
to $g_0$ at some $r_0\in[a,b]$ ; then equality holds if and only if $g=g_0$.
\end{theorem}

We note that our definition of mass is similar to (but different from) the definition of quasi-local
mass given by Brown and York \cite{YorkBrown}.
As a corollary we obtain some positive mass type or Penrose inequality type theorem.

\begin{corollary}
Under the hypothesis of theorem \ref{thintro}, if we also assume that
$$
\int_M\mathring{h}^{2}(R_{rr}-\mathring{R}_{rr})drd\mu_{\hat{G}}\geq 0,$$
then $\m(b)\geq \m(a)$, the equality case being the same as the one in
theorem \ref{thintro}.
\end{corollary}

This result is reminiscent of the positive mass theorem (see \cite{ShoenYau:masse} and
\cite{Witten:masse}), which states that under certain conditions (in particular  : the scalar
curvature belongs to $L^1(M)$) an asymptotically Euclidian manifold with nonnegative scalar
curvature has nonnegative ADM mass $m_{ADM}$. Moreover $m_{ADM}=0$ if and only if the
manifold is isometric to the Euclidian space (see also \cite{bartnik:qlm}, \cite{Wang:mass},
\cite{ChHerzlich}, \cite{ChNagy}, \cite{ChNagyATMP} for related results).

Of course our result is not really a positive mass theorem in dimension greater than
or equal to 3, because the usual energy condition  on the scalar curvature $R(g)\geq R(g_0)$
is replaced here by an integral condition on  the Ricci curvature. The similarity with the
Penrose inequality (see \cite{Penrose:inequality}) appears more clearly when $N(a)$ is a minimal
hypersurface of $g$,
in which case $\m(a)=(n-1)\mathring{h}(a)\mathring{h}'(a)\hat{\mu}(N)$, where $\hat{\mu}(N)$
is the volume of $N$ with respect to $\hat{G}$.

We also obtain some rigidity results in the spirit of those obtained in
\cite{Michel-Sarih}, \cite{Michel} or \cite{Sarhi} (see also \cite{Hopf}, \cite{Green}
or \cite{BurnsGerhard} for instance) : a complete Riemannian metric on $\R^n$ without
conjugate points and integrable Ricci curvature has nonpositive total scalar curvature
$\kappa$, and $\kappa$ vanishes if and only if the metric is flat.

\begin{corollary}
Under the hypothesis of theorem \ref{thintro}, and if $\m(b)\leq \m(a)$, then
$$
\int_M\mathring{h}^{2}(R_{rr}-\mathring{R}_{rr})drd\mu_{\hat{G}}\leq 0,$$
with the same equality case as in theorem \ref{thintro}.
\end{corollary}

Let us remark that when $R_{rr}\geq \mathring{R}_{rr}$  the Riccati equation
\eq{riccadni} shows that the map $r\mapsto m(r,x)$ is a non decreasing function
and gives an alternative proof of Bishop-Gromov type theorems (see the survey
\cite{Wei:Riccibound} for instance).

It is also interesting to note that when
$\int_N\mathring{h}^2(R_{rr}-\mathring{R}_{rr})d\mu_{\hat G}\geq0$
the Riccati equation \eq{riccadni} also shows that the map $r\mapsto \m(r)$ is
a non decreasing function. Moreover $\m$ will be constant if and only if
$g$ is as in the equality case of theorem \ref{thintro}.

We will first give a proof in the two-dimensional case in the next section. The
calculations are more straightforward and $g_0$ need not be a warped product. We
deal with the n-dimensional case in section 3. We then apply our result to
standard models. Finally, in the last section, we start a discussion about the
geometric invariance of our notion of mass.
\\

\noindent{\bf Acknowledgement} We are grateful to E. Aubry and Ph. Castillon for
interesting discussions about the n-dimensional case.

\section{The two-dimensional case}\label{lasection}

%

Let $(a,b)$ be a non-empty interval of $\R$.
We consider the surface $$M=(a,b)\times N,$$ $N$ being either
a closed interval of $\R$ or $\Sphere^1$.
Let  $g_0$ be a fixed background metric on $M$, with
$$
g_0=dr^2+h_0^2(r,\theta)d\theta^2,
$$
and let $g$ be another Riemannian metric on $M$, also of the form
\bel{g}
g=dr^2+h^2(r,\theta)d\theta^2.
\ee
We define two quantities
$$
m{(r,\theta)}=-h_0^2\partial_r (\ln(\frac{h}{h_0})),
$$
and
$$
\m(r)=\int_N m(r,\theta)d\theta.
$$

We shall call $\m(r)$  the {\it relative radial mass} of $g$ to $g_0$ at $r$.
Note that for $i \in \{ a,b \}$, $\m(i)$ does exist if we assume that near $i$,
$$
m(r,\theta)=m(\theta)+o(1),
$$
where $m$ is an integrable function on $N$ and the $o(1)$ is uniform with respect
to $\theta$.

Let us denote by $K$ and $K_0$  the Gau\ss $\mbox{ }$ curvature
of $g$ and $g_0$ respectively, so that
$$
K=-\frac{\partial^2_rh}{h},\;\; K_0=-\frac{\partial^2_rh_0}{h_0}
$$


\begin{lemma}\label{Riccati}
We  have  the  Riccati type formula
\bel{laformuleder}
\partial_r m=h^2_0 (K-K_0)+h_0^{-2}m^2.
\ee
\end{lemma}

\begin{proof}
We will denote by a $'$ the partial derivative relative to $r$. Let
$$v:=\frac{h'}{h}\;, \;\;v_0:=\frac{h'_0}{h_0}\;,\;\; V:=v-v_0=(\ln(\frac{h}{h_0}))'.$$
We have
$$V'=-K+K_0-v^2+v_0^2=-(K-K_0)-V^2+2v_0^2-2v_0v=-(K-K_0)-V^2-2v_0V,
$$
Derivating the relation  $m=-h_0^2V$ then gives
\begin{eqnarray*}
m'&=&-2h_0h'_0V-h_0^2V'=-2h_0^2v_0V-h_0^2V'\\
&=&(K-K_0)h_0^2+h_0^2V^2\\
&=& (K-K_0)h_0^2+m^2h_0^{-2}.
\end{eqnarray*}
\end{proof}

%

Recall that $g$ is asymptotic to $g_0$ at $i\in[a,b]$ if near $i$,
$h=h_0(1+o(1))$, the $o(1)$ being uniform relatively to $\theta$.

By integrating the formula \eq{laformuleder} on $M$, we thus obtain the

\begin{theorem}\label{rigidite}
If $(K-K_0)\in L^1(h^2_0drd\theta)$ and $\m$
has a limit at $a$ and $b$ then $(\frac{m}{h_0})^2\in L^1(drd\theta)$ and
$$
\m(b)-\m(a)\geq \int_{(a,b)\times N}(K-K_0)h_0^2drd\theta,
$$
with equality if and only if $K\equiv K_0$ and $m\equiv 0$ (so that $\frac{h}{h_0}$
does not depend on $r$). In particular, if $g$ is asymptotic to $g_0$ at some
$r\in[a,b]$ then we have equality if and only if $g=g_0$.
\end{theorem}

\begin{remark}
Some interesting results are also derived if we directly use the equation
\eq{laformuleder} or if we integrate it on $N$ only or on $(a,b)$ only.

As a first example let us assume that  $K(.,\theta)\geq K_0(.,\theta)$ on $(a,b)$ ; then
$m(.,\theta)$ is non-decreasing on $(a,b)$, and  $m(.,\theta)$ is constant if and
only if $K(.,\theta)= K_0(.,\theta)$ and $m(.,\theta)=0$ (so that
$\frac{h(.,\theta)}{h_0(.,\theta)}$ only depends on $\theta$).

Also note that under the assumption that $\int_N(K-K_0)h_0^2d\theta \geq 0$, we can conclude
that $\m$ is a non-decreasing function on $(a,b)$, and that $\m$ is constant if
and only if $K=K_0$ and $m=0$ (so that $\frac{h}{h_0}$ only depends on $\theta$).

Finally if $(K(.,\theta)-K_0(.,\theta))\in L^1(h^2_0(.,\theta)dr)$ and $m(.,\theta)$
has a limit at $a$ an $b$ then $(\frac{m}{h_0})^2(.,\theta)\in L^1(dr)$ and
$$
m(b,\theta)-m(a,\theta)\geq\int_{(a,b)}(K(.,\theta)-K_0(.,\theta))h_0^2dr,
$$
with equality if and only if $K(.,\theta)=K_0(.,\theta)$ and $m(.,\theta)=0$
(so that $\frac{h(.,\theta)}{h_0(.,\theta)}$ only depends on $\theta$).

\end{remark}

%
%
%
%
%


%
%
%

%

\begin{remark}
If $h_0$ depends only on $r$, let us define the function
$$
u(r,\theta)=h_0(r)\partial_r\ln(\frac{h}{h_0})(r,\theta),\;\;
$$
and  the one form  $U(r,\theta)=u(r,\theta)dr$.

A simple calculation shows that $m^2= h_0^2||U||^2_{g_0}$ and that
$$
\div_{g_0}U=-(\partial_ru+h_0^{-1}(\partial_rh_0) u)$$
so we get
$$
\div_{g_0}U=h_0(K-K_0)+h_0^{-1} ||U||^2_{g_0},
$$
which is a divergence form of the formula \eq{laformuleder}.
\end{remark}

\section{The $n$-dimensional case }

Let $M=(a,b)\times N$, where $N$ is a compact $(n-1)$-dimensional manifold, and
$(a,b)$ an open interval of $\R$. We will consider coordinate systems adapted
to  $M$, that is of type $(x^0,x^1,....,x^{n-1})=(x^i)=(r,x^A)$,
where $(x^1,...,x^{n-1})=(x^A)$ is a local coordinate system on $N$ (indices
with capital letters are used to denote coordinates on N).
On $M$, we consider a Riemannian  metric of the form
\bel{gdn}
g=dr^2+G(r)=dr^2+G_{AB}(r,x^D)dx^Adx^B,
\ee
where $G(r)$ is a family of Riemannian metrics on $N$.
We also consider a background warped product metric on M :
\bel{g0dn}
g_0=dr^2+\mathring{G}(r)=dr^2+\mathring{G}_{AB}(r,x^D)dx^Adx^B,
\ee
where $\mathring{G}(r)=\mathring{h}^2(r)\hat{G}$, $\mathring{h}$ being a positive
function on  $(a,b)$ and $\hat{G}$ a fixed metric on $N$.

We respectively denote by $S$ and $\mathring{S}$ the $g$ and $g_0$  Weingarten
maps (shape operators) of the level sets of $r$ with respect to the outward normal
$\partial_r$, whereas $II$ and $\mathring{II}$ are the associated second
fondamental forms.
In coordinates we have :
$$II_{AB}=\frac{1}{2}G'_{AB},
$$
$$
S_{A}^B=\frac{1}{2} G^{BC}G'_{AC},
$$
$$\mathring{II}_{AB}=\frac{1}{2}\mathring{G}'_{AB}=\mathring{h}\mathring{h}'\hat{G}_{AB},
$$
$$
\mathring{S}_{A}^B=\frac{1}{2} \mathring{G}^{BC}\mathring{G}'_{AC}=\mathring{h}^{-1}\mathring{h}'\delta^{B}_A,
$$
The Riemann curvature in two $\partial_r$ directions gives as usual some
Riccati equations for $S$ and $\mathring{S}$:
\bel{RiccaS}
R^B_{rAr}=-S'^B_A-S^C_AS^B_C,
\ee
\bel{RiccaS0}
\mathring{R}^B_{rAr}=-\mathring{S}'^B_A-\mathring{S}^C_A\mathring{S}^B_C,
\ee
and their traces
\bel{RiccaH}
R_{rr}=-H'-|S|^2_G,
\ee
\bel{RiccaH0}
\mathring{R}_{rr}=-\mathring{H}'-|\mathring{S}|^2_{\mathring{G}},
\ee
where $H=\Tr S$ and $\mathring{H}=\Tr \mathring{S}=(n-1)\mathring{h}^{-1}\mathring{h}'$
are the mean curvatures of the level sets of $r$ relative to $g$ and $g_0$ respectively.
As
$$|\mathring{S}|^2=(n-1)(\mathring{h}^{-1}\mathring{h}')^2=\frac{1}{n-1}\mathring{H}^2
\ {\rm and }\ |S|^2=|s|^2+\frac{1}{n-1}H^2,$$
where $s=S-\frac{1}{n-1}H\;I$ is the trace-free part of $S$, we finally obtain
\bel{RiccaHH0}
\mathring{H}'-H'= R_{rr} -\mathring{R}_{rr} + |{S}|^2 -|\mathring{S}|^2 =
R_{rr} - \mathring{R}_{rr} + \frac{1}{n-1}(H^2 - \mathring{H}^2) + |s|^2.
\ee
Recall that
$$
m(r,x^A)=\mathring{h}^{2}(\mathring{H}-H),
$$
and note that
$$(\mathring{h}^{2})'=2\mathring{h}\mathring{h}'=\frac{2}{n-1}\mathring{h}^{2}\mathring H.
$$
Combining this two relations with $\eq{RiccaHH0}$ gives
\bel{riccadn}
m'=\mathring{h}^{2}(R_{rr}-\mathring{R}_{rr})+\frac{1}{n-1}\frac{m^2}{\mathring{h}^2}
+\mathring{h}^{2}|s|^2
\ee

As in dimension two, we define $\m$ as
$$
\m(r)=\int_N\mathring{h}^{2}(\mathring{H}-H)d\mu_{\hat{G}}.
$$

In order to prove the inequality of the theorem \ref{thintro}, it suffices
to integrate the equation \eq{riccadn} over $M$.

We now deal with the equality case. We have
$$
\int_M(\frac{1}{n-1}\frac{m^2}{\mathring{h}^2}
+\mathring{h}^{2}|s|^2)drd\mu_{\hat{G}}=0,
$$
so that $m \equiv 0$ (thus $\mathring{H}=H$) and $s \equiv 0$.
Coming back to equation \eq{riccadn} we also have $R_{rr}=\mathring{R}_{rr}$.
Now $s=0$ and $H=\mathring H$ give
$$
S=\frac{1}{n-1}HI=\frac{1}{n-1}\mathring{H}I=\mathring{S}.
$$
Equations \eq{RiccaS} and \eq{RiccaS0} then imply that
$\mathring{R}^B_{rAr}={R}^B_{rAr}$.
Finally  $\mathring{S}={S}$ is equivalent to $(\mathring{h}^{-2}G_{AB})'=0$ so
that $G_{AB}=\mathring{h}^{2}\w{G}_{AB}$ for some metric $\w{G}$ on $N$. If $g$
is asymptotic to $g_0$ at some $r_0\in[a,b]$ (i.e.
$||G-\mathring{G}||_{\mathring{G}}=o(1)$ near $r_0$) then $\w G=\hat G$ so that
$g=g_0$ and this concludes the proof of theorem \ref{thintro}.

\begin{remark} \label{remdet}

We denote by $$h=\sqrt{\det g}=\sqrt{\det G}$$ and  by
$$h_0=\sqrt{\det g_0}=\sqrt{\det \mathring{G}}=\mathring{h}^{n-1}\sqrt{\det \hat{G}}$$
the volume element of $g$ and $g_0$ respectively. A standard calculation (making use
of the formula for the derivative of a determinant) shows that
$$h'=\frac{1}{2}h\Tr_g(g')=\frac{1}{2}h\Tr_G(G')=hH$$
and  $h_0'=h_0\mathring{H}$, so we have
\bel{mdet}
m=\mathring{h}^2(\ln(\frac{h_0}{h}))',
 \ee
which is the definition of $m$ in the two-dimensional case because $\mathring{h}=h_0$ in
that case.
\end{remark}

\begin{remark}\label{remdiv} As in dimension 2, the formula \eq{riccadn} can be written
in divergence form. To do this, consider the fonction $u=\mathring{h}^{3-n}(H-\mathring{H})$
and the one-form $U=u\;dr$. The divergence of $U$
is
\begin{eqnarray*}
\div_{g_0}U&=&-u'-\mathring H u\\
&=&(3-n)\mathring{h}^{2-n}\mathring{h}'(\mathring{H}-H)+\mathring{h}^{3-n}(\mathring{H}-H)'
+\mathring{h}^{3-n}\mathring H (\mathring{H}-H)\\
&=&\mathring{h}^{3-n}(\mathring{H}-H)'+\frac{2}{n-1}\mathring H\mathring{h}^{3-n}(\mathring{H}-H)\\
&=&\mathring{h}^{3-n}(R_{rr}-\mathring{R}_{rr})+\frac{1}{n-1}{||U||^2_{g_0}}{\mathring{h}^{n-3}}
+\mathring{h}^{3-n}|s|^2,
\end{eqnarray*}
the last equality coming from \eq{RiccaHH0}.
\end{remark}

\section{Some simple conditions for the relative radial mass to be
well defined}

In this section, we give sufficient conditions for the radial mass $\m(i)$ to
be well defined when $i \in \{ a,b \}$. This will be the case if $g$ is in
some way asymptotic  to $g_0$ at $i$. The most simple case is certainly when
near $i$,
\bel{scm}
m(r,x)=m(x)+o(1),
\ee
$m$ being a continuous function on $N$ ; then
$$
\m(i)=\int_Nm(x)d\mu_{\hat{G}}(x).
$$

The aim of this section is to give some simple asymptotic conditions on
$g$ and $g_0$ which imply the convergence of $\m$. We recall the formula
$$
m(r,x)=-\frac12\mathring{h}^2\Tr(G^{-1}G'-\mathring{G}^{-1}\mathring{G}').
$$

We first give some conditions which are natural when for instance one wants
to define a mass at infinity and $\mathring{h}=r$ or $\mathring{h}=e^r$. Let
us assume that $\int^i\mathring{h}^{-2}(r)dr$ converges and consider the
function
$$
F_i(r)=\int^i_r\mathring{h}^{-2}(r)dr=o(1).
$$
Assume that near $i$,
 \bel{sc1}
G(r)=\mathring{h}^2[\hat{G}+F_i\w G],
\ee
where $\w G$   and $\mathring{h}^2F_i(\w{G})'$ are bounded. In this case
we obtain that near $i$,
$$
\mathring{h}^2(G^{-1}G'-\mathring{G}^{-1}\mathring{G}')=\hat{G}^{-1}[\w G+\mathring{h}^2F_i(\w{G})']+O(F_i)+O(F_i^2\mathring{h}\mathring{h}'),
$$
so that if $F_i^2\mathring{h}\mathring{h}'=o(1)$, then
$$
m(r,x)=-\frac12\Tr_{\hat{G}}[\w G+\mathring{h}^2F_i(\w{G})']+o(1).
$$
In particular, if $\w G+\mathring{h}^2F_i(\w{G})'$ is continuous up to $i$, then
$$
\m(i)=-\frac12\int_{N}\Tr_{\hat{G}}[\w G+\mathring{h}^2F_i(\w{G})'](i,.)d\mu_{\hat{G}}.
$$

We now give conditions which are natural if we want to define a mass at $r=0$,
with $g_0$ written in polar coordinates and $g$ possibly singular at the base
point. If $\int^i\mathring{h}^{-2}(r)dr$ diverges and if we assume that near $i$,
\bel{sc1b}
G(r)=\mathring{h}^2[\hat{G}+\overline{G}],
\ee
with $\overline{G}=o(1)$, $\mathring{h}\mathring{h}'\overline{G}=o(1)$, and
$\mathring{h}^2\overline{G}'=O(1)$,
then
$$
\mathring{h}^2(G^{-1}G'-\mathring{G}^{-1}\mathring{G}')
=\mathring{h}^2\hat{G}^{-1}\overline{G}'+o(1).
$$
In particular if $\mathring{h}^2\hat{G}^{-1}\overline{G}'$ is continuous up to
$i$ then
$$
\m(i)=-\frac12\int_{N}\mathring{h}^2\Tr_{\hat{G}}\overline{G}'(i,.)d\mu_{\hat{G}}.
$$

\section{Applications to the rigidity of some models}

In this section we apply theorem \ref{thintro} to the case where the Riemannian
manifold $(M,g_0)$ is a model space of constant radial Ricci curvature. The equality case in
theorem \ref{thintro} then provides rigidity results for the model $g_0$ if the
metric $g$ is asymptotic to $g_0$ at some $r\in[a,b]$. Note that our results differ
from the usual positive mass theorems, as we do not require $b=+\infty$, nor do
we need the two metrics to be asymptotic at $b$.

Recall that we assume that
$$
g_0=dr^2+\mathring{h}^2 \hat{G} \ {\rm and} \ g=dr^2+G(r),
$$
so that in particular
$$
\mathring R_{rr}=-(n-1)\frac{\mathring{h}''}{\mathring h}.
$$

\subsection{Hyperbolic type Metric}

We consider the case where $\mathring{h}=e^r+ke^{-r}$ and $M=(a,+\infty)\times N$
with $k\in \R$ and $a$ is the zero of $\mathring{h}=e^r+ke^{-r}$ (or $a=-\infty$
if  $\mathring{h}$ never vanishes). By an immediate computation,
$\mathring R_{rr}=-(n-1)$, and theorem \ref{thintro} becomes the following

\begin{corollary}\label{hyperbolique}

If the limit of $\m$ exists at $a$ and $+\infty$
and $\mathring{h}^{2}(R_{rr}+(n-1))\in L^1(drd\mu_{\hat{G}})$  then
$$
\m(+\infty)-\m(a)\geq \int_M\mathring{h}^{2}(R_{rr}+(n-1))drd\mu_{\hat{G}}$$
with equality if and only if $g$ is a warped product ($G=\mathring{h}^2\w{G}$,
$\w{G}$ being a fixed metric on $N$). Moreover if $g$ is asymptotic to $g_0$
at some $r_0\in[a,+\infty]$, then equality  occurs if and only if $g=g_0$.
\end{corollary}

\subsection{Euclidian type Metric}

In this case, $\mathring{h}=r$ and $M=(0,+\infty)\times N$. As $\mathring R_{rr}=0$
we obtain the

\begin{corollary}\label{euclidien} If the limit of $\m$ exists at $0$ and $+\infty$ and
$r^{2}R_{rr}\in L^1(drd\mu_{\hat{G}})$ then
$$
\m(+\infty)-\m(0)\geq \int_Mr^{2}R_{rr}drd\mu_{\hat{G}}$$
with equality if and only if $g$ is a warped product ($G=r^2\w{G}$, $\w{G}$ being
a fixed metric on $N$). Moreover if $g$ is asymptotic to $g_0$ at some $r_0\in[0,+\infty]$
then equality  occurs if and only if $g=g_0$.
\end{corollary}

\subsection{Cylindrical type Metric}

This is the case where $\mathring{h}=1$ and $M=(-\infty,+\infty)\times N$. Once again,
$\mathring R_{rr}=0$ and so we have the

\begin{corollary}\label{cylindre}
If the limit of $\m$ exists at $-\infty$ and $+\infty$ and $R_{rr}\in L^1(drd\mu_{\hat{G}})$
then
$$
\m(+\infty)-\m(-\infty)\geq \int_M R_{rr} drd\mu_{\hat{G}},
$$
with equality if and only if $g$ is a  product. Moreover if $g$ is asymptotic to $g_0$
at some $r_0\in[-\infty,+\infty]$ then equality occurs if and only if $g=g_0$.
\end{corollary}

\subsection{Spherical type Metric}

In this final case, $\mathring{h}=\sin (r)$ and $M=(0,\pi)\times N$, so that
$\mathring R_{rr}=(n-1)$ and theorem \ref{thintro} yields

\begin{corollary}\label{Sphere} If the limit of $\m$ exists at $0$ and $\pi$ and
$\sin^2r(R_{rr}-(n-1))\in L^1(drd\mu_{\hat{G}})$  then
$$
\m(\pi)-\m(0)\geq \int_M (R_{rr}-(n-1))\sin^2r drd\mu_{\hat{G}}$$
with equality if and only if $g$ is a  warped product (with $G=\sin^2(r)\w{G}$,
$\w{G}$ being a fixed metric on $N$). Moreover if  $g$ is asymptotic to $g_0$ at
some $r_0\in[0,\pi]$ then equality occurs if and only if $g=g_0$.
\end{corollary}

\section{About the geometric invariance of the mass}

As it is usual when a definition of mass is given in a particular system of coordinates,
we have to check that (maybe under certain conditions on $g$ and $g_0$) this notion does not
depend on the coordinate system chosen to define it. This problem can be formulated in two
different ways. The first one, which is more natural in our context, will be treated here
and the answer turns out to be positive when $n=3$. We did not investigate any further
the second one, which is more natural when dealing with positive mass theorems. At this point
let us simply note that the method used in \cite{ChHerzlich} and \cite{ChNagy} combined
with a divergence formula (like the one in remark \ref{remdiv}) might be useful to solve
this problem.

The first way to understand the problem is as follows. Let $\mathcal X$ be a manifold endowed
with two Riemanniann metrics $g$ and $g_0$ and let $E$ be one end of $\mathcal X$. We assume
that there exists a coordinate system $(r,x)$ on $E$ making it diffeomorphic to
$M:=[a,+\infty)\times N$, and in which $g$ and $g_0$ take the form \eq{gdn} and \eq{g0dn}
respectively. Let $\m(+\infty)$ be the relative mass at infinity of $g(r,x)$ relative to
$g_0(r,x)$. Now suppose that there exists another coordinate system $(\w{r},\w{x})$ of $E$
making it diffeomorphic to $\w M:=[\w a,+\infty)\times N$, and in which $g$ takes the form
$$
\w{g}:=d\w{r}^2+\w{G}(\w r),
$$
while $g_0$ takes the form
$$
\w{g}_0:=d\w{r}^2+\w{\mathring{G}}(\w r),
$$
where $\w{\mathring{G}}(\w r)=\w{\mathring{h}}(\w r)^2\w{\hat{G}}$ and $\w{\hat{G}}$
is a fixed metric on $N$. Let $\w{\m}(+\infty)$ be the mass at infinity of
$\w{g}(\w r,\w x)$ relative to $\w g_0(\w r,\w x)$. The problem in to find some natural
conditions on $g_0$ and $g$ which ensure that $\m(+\infty)=\w{\m}(+\infty)$.

Let us denote by $\Phi$ the change  of coordinates :
$$(r,x)=\Phi(\w r,\w x),$$
so that $\w g=\Phi^* g$ and $\w g_0=\Phi^* g_0$. We assume that $\Phi$ is orientation preserving.

\begin{lemma}\label{rinfini} Near infinity, the first projection of the change of
coordinates $\Phi$ satisfies
$$|r-\w r|\leq  C,$$
where $C$ is a positive constant. A straightforward consequence is that $r\rightarrow+\infty$
if and only if $\w r\rightarrow+\infty$.
\end{lemma}

\begin{proof}
The form of $g_0$ and $\w{g}_0$ and the fact that $\Phi$ and its inverse are isometries
between the two metrics imply that $\partial_r{\w{r}}\leq 1$ and also $\partial_{\w{r}}r\leq 1$.
So for all $\w r\geq \w r_0$ we have
$$
r(\w r,\w{\theta})-r(\w r_0,\w{\theta})= \int_{\w r_0}^{\w r}\partial_{\w r}r(s,\w{\theta})ds
\leq (\w r-\w r_0),
$$
thus
$$
r\leq \w r+ {C}(\w{\theta}),
$$
where ${C}(\w{\theta})$ is continuous on $N$. Now as $N$ is compact, ${C}(\w{\theta})$ is bounded
and we have
$$
r\leq \w r+ C,
$$
where $C$ is a positive constant. In the same way, increasing $C$ if necessary,
$$
\w r\leq r+ C,
$$
\end{proof}

Set $h=\sqrt{\det g}$ and  $h_0=\sqrt{\det g_0}$ and consider the function
$$
v=\ln(\frac{h_0}{h}).
$$
Note that $v$ does not depend on any coordinate system, whereas $h$ and $h_0$ do. Let us denote by
$\Delta_0$ the Laplacian  relative to $g_0$.

\begin{proposition}\label{mgeom} If $N$ is without boundary, $n=3$ and $\Delta_0 v\in L^1(E, g_0)$
then $\m(+\infty)=\w{\m}(+\infty)$.
\end{proposition}

\begin{proof}
Consider the one-form $$\omega=dv.$$
We have $\div_0\omega=\Delta_0v$ and if we denote by $\nu$ (resp. $\w{\nu}$) the outgoing
$g_0$-normal to $N(r)$ (resp. $\w N(\w r )$) , we have
$\langle \omega,\nu\rangle_{g_0}=\partial_r v$
(resp. $\langle \omega,\w{\nu}\rangle_{g_0}=\partial_{\w r} v$).
Let $\Omega_{r,\w r}$ be the open relatively compact set with boundary $N(r)\cup \w N(\w r)$.
As a consequence of the divergence theorem, we have
$$
\int_{N(r)}(\partial_r v)\mathring{h}^{n-1}d\mu_{\hat{G}}
- \int_{\w N(\w r)}(\partial_{\w r} v)\w{\mathring{h}}^{n-1}d\mu_{\w{\hat{G}}}
=\int_{\Omega_{r,\w r}} \epsilon \div_0 \omega \;d\mu_{g_0},
$$
where $\epsilon = \pm 1$ is constant on every connected component of $\Omega_{r,\w r}$.

Now from the lemma \ref{rinfini}, when $r\rightarrow+\infty$ the same property holds for
$\w r$, so that the set $\Omega_{r,\w r}$ goes to infinity. But as $\div_0 \omega\in L^1(E)$
the last integral goes to zero, and we obtain
$$
\m(r)=\w{\m}(\w r)+o(1).
$$
\end{proof}

\begin{remark}
If $\m(r)$ converges at infinity and $R_{rr}-\mathring{R}_{rr}\in L^1(E)$ (as in theorem
\ref{thintro} with $n=3$) then $\Delta_0 v \in L^1(E)$, thus the $L^1$ condition in
proposition \ref{mgeom} is natural. In order to justify this affirmation, we compute :
\begin{eqnarray*}
\Delta_0 v&=&-h_0^{-1}\partial_i(g_0^{ij}h_0\partial_j v)\\
&=&-h_0^{-1}\partial_r(h_0\partial_r v)-h_0^{-1}\partial_A(\mathring{G}^{AB}h_0\partial_B v)\\
&=&-{\mathring{h}}^{1-n}\partial_r({\mathring{h}}^{n-1}\partial_r v)
+{\mathring{h}}^{-2}\hat{\Delta}v,\\
\end{eqnarray*}
where $\hat{\Delta}$ is the Laplacian relative to $\hat{G}$.
Now the fact that $\int_N\hat{\Delta}vd\mu_{\hat{g}}=0$ and that $m=\mathring{h}^{n-1}\partial_r v$
when $n=3$ combined with the formula \eq{riccadni} concludes the proof.
\end{remark}

\begin{remark}
The change of variable $(r,\theta)=(\w r,\lambda \w{\theta})$  gives
$\w{\mathring{h}}=\lambda\mathring{h}$ and $\m(+\infty)=\lambda^{3-n}\w{\m}(+\infty)$ ; this shows
that when $n\neq 3$, our mass is not a geometric invariant.
\end{remark}

Let us now formulate another problem also related to the geometric invariance of the mass. Let
$(\mathcal X,g)$ be a Riemannian manifold. Let $E$ be one end of $X$. We assume that there exists
coordinate system $(r,x)$ of $E$ making it diffeomorphic to $M:=(a,+\infty)\times N$, and in which
$g$ takes the form \eq{gdn}. Let $\m(+\infty)$ be the mass of $g(r,x)$ relative to $g_0(r,x)$.
Take another coordinate system $(\w{r},\w{x})$ of $E$ making it diffeomorphic to
$\w M:=(\w a,+\infty)\times N$, and in which $g$ takes the form
$$
\w{g}:=d\w{r}^2+\w{G}(\w r).
$$
Let $\w{\m}(+\infty)$ be the mass at infinity of $\w{g}(\w r,\w x)$ relative to $\w g_0(\w r,\w x)$.
Once again, we ask for natural conditions on $g_0$ and $g$ that guarantee that $\m(+\infty)=\w{\m}(+\infty)$.
If they exist, $\m(+\infty)$ is, in this sense, a geometric invariant.

Let $\Phi$ be the transition of coordinates, so $\w g=\Phi^*g $. The fundamental difference
with the previous case is that we do not have $\w g_0=\Phi^* g_0$ in general. As in \cite{ChNagy}
or \cite{ChHerzlich}, we may hope that if $g_0$ is well chosen and $g$ asymptotic to $g_0$, we will
be able to investigate the asymptotical behaviour of $\Phi$ and show that the mass at infinity is a
geometric invariant.

We will only deal with two very particular cases where we can prove that our mass corresponds
to the usual one (and is therefore geometric). In section \ref{CH} we establish in dimension 2, in
the asymptotically hyperbolic context, our mass at infinity is the time component of the energy
momentum vector of Chrus\'ciel and Herzlich (see also Wang \cite{Wang:mass}), so that its sign
is a geometric invariant. We then show that in dimension 3, in the asymptotically flat context,
 our mass at infinity is (proportionnal to) the usual ADM mass (see section \ref{adm}).

\subsection{Comparison with the mass of asymptotically hyperbolic surfaces}\label{CH}

In the two dimensional asymptotically hyperbolic context, so when $h_0=e^r+ k e^{-r}$ (k=-1,0,1), we will
see that our mass $\m(+\infty)$ at infinity correspond to the time component
$p_{(0)}$ of the energy momentum vector of  Chrus\'ciel and Herzlich \cite{ChHerzlich}
up to a positive multiplicative constant.
To see that, we will  compute  the vector field $\U=\U(V)$ they use with $V=V_{(0)}$.
$$
\U^i=2\sqrt{\det g}(Vg^{i[k}g^{j]l}\mathring{D}_je_{kl}+D^{[i}Vg^{j]k}e_{jk}).
$$
In our context, $g=dr^2+h^2d\theta^2$,  $g_0=dr^2+h_0^2d\theta^2$,
$e=g-g_0=(h^2-h_0^2)d\theta^2$, so
$$
\U^i=2\sqrt{\det g}(Vg^{i[\theta}g^{j]\theta}\mathring{D}_je_{\theta\theta}
+D^{[i}Vg^{\theta]\theta}e_{\theta\theta}).
$$
First $2D^{[i}Vg^{\theta]\theta}=D^iVg^{\theta\theta}-D^1Vg^{i\theta}=h^{-2}\delta^i_rD^rV$
and the quantity
$2[g^{i[\theta}g^{j]\theta}]=g^{i\theta}g^{j\theta}-g^{ij}g^{\theta\theta}$ is non
trivial iff $i=j=r$ and then its value is $-h^2$. So $\U^\theta=0$ and
\begin{eqnarray*}
\U^r&=&h^{-1}(e_{\theta\theta}D^rV-V\mathring D_re_{\theta\theta})
=h^{-1}[e_{\theta\theta}\partial_rV-V(\partial_re_{\theta\theta}-2h_0h'_0e_{\theta\theta})]\\
&=&h^{-1}[(h^2-h_0^2)\partial_rV-2V(h\partial_rh-h^2h_0^{-1}h_0')]
\end{eqnarray*}
Now we know that $V=h_0+O(e^{-r})$ and $V'=h_0'+O(e^{-r})=h_0+O(e^{-r})$. Assume that we have
$h=h_0+e^{-r}\cH$, where $\cH$ is bounded and also $\partial_rh=h_0'+e^{-r}(\partial_r\cH-\cH)$
where $\partial_r\cH-\cH$ is bounded, so that $g$ is in some sense asymptotical to $g_0$.
Then up to terms that cancel at $r=+\infty$, we have
$$
\U^r=h^{-1}(h^2-h_0^2)h_0'-2h_0^2\partial_r(\frac{h}{h_0})+o(1)=h^{-1}(h^2-h_0^2)h_0'+2m+o(1),
$$
where $m$ is our $m(r,\theta)$. Thus we have $m(r,\theta)=2\cH-\partial_r\cH+o(1)$ and
$h^{-1}(h^2-h_0^2)h_0'=2\cH+o(1)$.
Let
$$\w{m}(r)=\int_N \cH d\theta,$$
then
$$
\m(+\infty)=\lim_{+\infty}(2\w{m}-\w{m}').
$$
It is an easy exercice to show that in that case the only
possibility is $\lim_{+\infty}(\w{m}')=0$ and then
$$
\m(+\infty)=\lim_{+\infty}(2\w{m}).
$$
Finally we have obtained
$$
p_{(0)}=\lim_{+\infty}\int_N\langle U,\partial_r\rangle_{g_0}d\theta=\lim_{+\infty}\int_N\U^rd\theta=2\m(+\infty)+\m(+\infty)=3\m(+\infty).
$$

\subsection{Comparison with the ADM mass}\label{adm}
In dimension 3, in the asymptotically flat context, we will see
that our mass $\m(+\infty)$ is the ADM mass up to a positive multiplicative
constant.
In that context, we have $\mathring{h}=r$ and $\hat{G}$ is the round metric
on the two sphere $N$.
Recall first that the Hawking mass on $N(r)$ is equal to
$$
m_H(r)=\frac{A^{1/2}}{(16\pi)^{3/2}}\left(16\pi-\int_{N(r)}H^2d\nu_r\right),
$$
where  $d\nu_r$ is the $g$-induced mesure on $N(r)$, and $A$ is the
$g$-area of $N(r)$.
Recall also that when $r$ goes to infinity, then under asymptotically Euclidian
conditions, $m_H(r)$  goes to the ADM mass $m_{ADM}$ (see \cite{HuiskenIlmanen} for instance).
From our definition of $m(r,x)$, we may assume that near infinity
$$
H(r,x)=\mathring{H}(r,x)-\frac{m(r,x)}{r^2}+o(\frac{1}{r^2}).
$$
So we obtain
$$
\m(r)=8\pi m_H(r)+o(1)=8\pi m_{ADM}+o(1),
$$
which gives the desired result in the limit at infinity.

\bibliographystyle{amsplain}

\bibliography{../references/newbiblio,%
../references/reffile,%
../references/bibl,%
../references/hip_bib,%
../references/newbib,%
../references/PDE,%
../references/netbiblio,%
../references/erwbiblio,%
 stationary}

\def\cprime{$'$}
\providecommand{\bysame}{\leavevmode\hbox to3em{\hrulefill}\thinspace}
\providecommand{\MR}{\relax\ifhmode\unskip\space\fi MR }
\providecommand{\MRhref}[2]{%
  \href{http://www.ams.org/mathscinet-getitem?mr=#1}{#2}
}
\providecommand{\href}[2]{#2}
\begin{thebibliography}{10}

\bibitem{bartnik:qlm}
R.~Bartnik, \emph{New definition of quasilocal mass}, Phys.\ Rev.\ Lett.
  \textbf{62} (1989), 2346--2348.

\bibitem{YorkBrown}
J.D. Brown and J.W. {York, Jr.}, \emph{Quasilocal energy and conserved charges
  derived from the gravitational action}, Phys.\ Rev. \textbf{D47} (1993),
  1407--1419.

\bibitem{BurnsGerhard}
K.~Burns and G.~Knieper, \emph{Rigidity of surfaces with no conjugate points.},
  J. Differential Geom. \textbf{34} (1991), no.~3, 623--650.

\bibitem{ChHerzlich}
P.T. Chru\'{s}ciel and M.~Herzlich, \emph{The mass of asymptotically hyperbolic
  {R}iemannian manifolds}, Pacific Jour.\ Math. \textbf{212} (2003), 231--264,
  dg-ga/0110035.

\bibitem{ChNagy}
P.T. Chru\'{s}ciel and G.~Nagy, \emph{The {H}amiltonian mass of asymptotically
  {anti-de Sitter} space-times}, Class.\ Quantum Grav. \textbf{18} (2001),
  L61--L68, hep-th/0011270.

\bibitem{ChNagyATMP}
\bysame, \emph{The mass of spacelike hypersurfaces in asymptotically {anti-de
  Sitter} space-times}, Adv.\ Theor.\ Math.\ Phys. \textbf{5} (2002), 697--754,
  gr-qc/0110014.

\bibitem{Green}
L.W. Green, \emph{A theorem of E. Hopf}, Michigan Math. J. (1958), 31--34.

\bibitem{Hopf}
E.~Hopf, \emph{Closed surfaces without conjugate points.}, Proc. Nat. Acad.
  Sci. (1948), 47--51.

\bibitem{HuiskenIlmanen}
G.~Huisken and T.~Ilmanen, \emph{The inverse mean curvature flow and the
  riemannian penrose inequality.}, J. Differential Geom. \textbf{59} (2001),
  353--437.

\bibitem{Michel}
R.~Michel, \emph{Sur quelques probl\`emes de g\'eometrie globale des
  g\'eod\'esiques}, Bol. Soc. Bras. Mat. Vol. \textbf{9} (1978), no.~2, 19--38.

\bibitem{Michel-Sarih}
R.~Michel and M.~Sarih, \emph{La rigidit\'e des plans sans points conjugu\'es
  asymptotiquement euclidiens}, Arch. Math. \textbf{52} (1989), 500--506.

\bibitem{Penrose:inequality}
R.~Penrose, \emph{Naked singularities}, Ann.\ New York Acad.\ Sci. \textbf{224}
  (1973), 125--134.

\bibitem{Sarhi}
M.~Sarhi, \emph{Absence de points conjugués et courbure intégrale dans
  {R}$^n$}, Algebras Groups Geom. \textbf{22} (2005), no.~1, 95--107.

\bibitem{ShoenYau:masse}
R.~Shoen and S.~T. Yau, \emph{The proof of the positive mass theorem. {II}.},
  Commun. Math. Phys. \textbf{9} (1981), 231--260.

\bibitem{Wang:mass}
X.~Wang, \emph{The mass of asymptotically hyperbolic manifolds}, J.
  Differential Geometry \textbf{57} (2001), 273--299.

\bibitem{Wei:Riccibound}
G.~Wei, \emph{Manifolds with a lower ricci curvature bound},
  arXiv:math/0612107v1 (2006).

\bibitem{Witten:masse}
E.~Witten, \emph{A new proof of the positive energy theorem.}, Commun. Math.
  Phys. \textbf{80} (1981), 381--402.

\end{thebibliography}

\end{document}